\documentclass[12pt]{amsart}

\newcommand{\beq}{\begin{equation}}
\newcommand{\eeq}{\end{equation}}

\newtheorem{prop}{Proposition}

\newtheorem{rem}{Remark}

\usepackage{amsmath}

\title{{\bf Closed-Form Solution of Polynomial Equations   }} 
\begin{document}
\author{Alexander Kheyfits }

\markright{this page style empty}

\begin{abstract} Complex analysis is used to find closed-form expressions of the roots of a univariate polynomial of any degree as integrals of elementary functions.  \end{abstract} 

\maketitle

\section{Introduction} 
The roots of any univariate algebraic equation of degree four or less can be expressed explicitly through its coefficients as superpositions of four arithmetic operations and radicals, while for no $n\geq 5$ there exists a general formula that represents all the $n$ complex roots of every $n-$th degree polynomial through its coefficients by means of the arithmetic operations and radicals. The reader can consult, e.g., \cite{Di} or \cite{Ki} for the history and the details of these results. The topic is of interest even today; see, e.g., a comprehensive current monograph of Teo Mora \cite{Mo} and the very recent papers by Boyd \cite{Bo} and Nash \cite{Na}. 

The Abel-Galois result did not, of course, stop scientists from search for other means to represent the roots of the quintics and more general equations. Already in 1858, Brioschi, Hermite, and Kronecker independently expressed the roots of a general quintic through elliptic modular functions; quarter century later Klein published his famous book \cite{Kl}. Thus, in solving the quintics, mathematicians have followed the standard path of extending the existing toolbox. For instance, the equation $x^2=-1$ cannot be solved in real numbers, but it can be done if we extend the available numbers to the imaginary numbers; to study the ordinary differential equation $x^2y''+xy'+(x^2-\nu^2)y=0$, which cannot be solved in elementary functions, the special (Bessel) functions were introduced and included in the set of possible solutions. These examples go on and on. Indeed, a sextic equation can be solved in terms \\

\footnoterule   

\emph{Key words:} Univariate polynomials; Explicit formulas expressing the roots as integrals of elementary functions.  \\

\emph{2010 Mathematics Subject Classification}: 12D10; 12E12; 97H30.  

\pagebreak
\noindent of Kamp$\acute{e}$ de F$\acute{e}$riet functions, and a septic in terms of the hyperelliptic functions and associated theta functions of degree $3$. To the best of the author's knowledge, no information regarding algebraic equations of higher degree is currently available.  

In this note we solve any algebraic equation in quadratures, i.e., we find explicit representations of the roots of univariate polynomials of any degree as integrals of elementary functions along the positive real axis. To derive these integral formulas, we apply the method from \cite{AnIo, Kh, Kh1}. \\

\section{Integral Representations of Roots of Polynomials with Simple Roots}
It is convenient to us to represent a generic polynomial with any real or complex coefficients as 
\[P(z)=z^n+a_1z^{n-1}+a_2z^{n-2}+\cdots +a_{n-2}z^2+a_{n-1}z+a_n.\] 
If $P(z)$ has at least one pair of multiple roots, then the discriminant of $P$ is zero, and the multiplicities of all the roots can be determined algebraically, e.g., by Horner's method. Now, by making use of the Euclidean algorithm, we can reduce $P$ to a polynomial of smaller degree without multiple roots. Thus from now on, we suppose that the polynomials under consideration have only simple roots.

Denote the roots of $P(z)=0$ as $\xi_i, \, i=1,2,\ldots,n$, where all the $n$ roots $\xi_i$ are different complex (or maybe real) numbers. We can also assume that all $\xi_i \neq 0$, since otherwise the problem is again reduced to a polynomial of lesser degree. Set $M=1+ \max \{|a_j|,\, 1\leq j\leq n\}, $ and let $A$ be any number such that $A>M$. Introduce another polynomial
\[f(z)= P(z+A)= z^n+b_1z^{n-1}+b_2z^{n-2}+\cdots +b_{n-1}z+b_n .\]
The equations $P(z)=0$ and $f(z)=0$ are equivalent, however, all the roots of $f(z)=0$ have negative real parts; the goal of the shift $z\rightarrow z+A$ is to clear the positive semi-axis $x=\Re z \geq 0$ from the zeros of $f$. In particular, $b_n=f(0)\neq 0$. We also mention that if $A$ is large enough, then $b_n \approx A^n (1+a_1/A)$. \\

The equation $f(z)=0$ can be rewritten as
\beq z= \sqrt[n]{-(b_1z^{n-1}+b_2z^{n-2}+\cdots +b_{n-1}z+b_n)}.\eeq
To work witth single-valued holomorphic functions, we cut the complex $z-$plane along the positive 
$x-$axis, and in the slit plane consider a closed contour $\Omega$, consisting of two horizontal segments 
\beq z=x \pm \imath \epsilon,\, 0\leq x<R ,\eeq 
whose left ends are connected with a small semi-circle $|z|=\epsilon,\, \Re z<0$, and their right end-points $R \pm \imath \epsilon $ are connected with a large arc of the circle $|z|=R'$; if $R$ is big enough, the contour  contains all the roots of the polynomial $f(z)$.  \\

In the domain $D$ bounded by the contour $\Omega=\partial D$, all arguments satisfy 
$0< \arg w <2\pi$. Consider $n$ holomorphic functions
\[f_k(z)=z - \: _{(k)}\sqrt[n]{-(b_1z^{n-1}+b_2z^{n-2}+\cdots +b_{n-1}z+b_n)}  \]
$k=0,1,2,\ldots ,n-1, \; z\in D$, where $_{(k)}\sqrt[n]{} \;$ is the $k-th$ branch of the radical in  (1), i.e.,
\[ _{(k)} \sqrt[n]{\cdot}=\big|b_1z^{n-1}+b_2z^{n-2}+\cdots +b_{n-1}z+b_n \big|^{1/n} \vspace{.2cm}\times \]
\beq \times \exp\left\{-\frac{\imath}{n}\arg\left(b_1z^{n-1}+b_2z^{n-2}+\cdots +b_{n-1}z+b_n\right) \right\} 
e^{\frac{2k\pi}{n} \imath} .\eeq
In particular,
\[f_0(z)=z - |b_1z^{n-1}+b_2z^{n-2}+\cdots +b_{n-1}z+b_n|^{1/n} \vspace{.2cm}\times \]
\[ \times \exp\left\{-\frac{\imath}{n}\arg\left(b_1z^{n-1}+b_2z^{n-2}+\cdots +b_{n-1}z+b_n\right) \right\}. \vspace{.2cm}\]

We also consider a linear function 
\[h(z)=\alpha z+ \beta  \]
with positive coefficients $\alpha >0$ and $\beta >0$, which has one root in the left half-plane.

We will show that the coefficients $\alpha $ and $\beta $ can be chosen so that $|f_k(z)<|h(z)|$ on the contour 
$\Omega$ for every $k=0,1,\ldots, n-1$, hence we can apply Rouch\'e's theorem to each pair  $h$ and $f_k$ for every $k=0,1,\ldots, n-1$, and conclude that every function $f_k(z)$ has the unique root in the plane. \\

The contour $\Omega$ consists of three parts; {\bf a)} small semi-circle of radius $\epsilon$ in the left half-plane; {\bf b)} two parallel horizontal segments $z=x\pm \imath \epsilon,\; 0<x<R'$; and {\bf c)} a big arc of the circle $|z|=R$, connecting the right ends of these horizontal segments. \\

{\bf a)} If $|z|=\epsilon, \; \epsilon>0$ and $\Re z<0$,  or $0\leq \Re z\leq 1$, then 
\[|f_k(z)| \leq \epsilon +\sqrt[n]{2B},\]
where $B=\max |b_j|,\, 1\leq j\leq n$. On the other hand, here $|h(z)|\geq \beta - \epsilon \alpha$, hence we must choose the parameters to satisfy the inequality
\beq \sqrt[n]{2B} +(\alpha +1)\epsilon \leq \beta. \eeq

{\bf b)} In this case $|f_k(z)|\leq \sqrt[n]{nB}x$, while 
\[|h(z)| \geq \alpha(x-\epsilon)+\beta \geq \sqrt[n]{nB}x. \]
Hence, in this case the parameters have to satisfy two inequalities $\alpha\geq \sqrt[n]{nB}$ and 
$\beta\geq \alpha \epsilon$.  \\

{\bf c)} Finally, if $|z|= R$, then $|f_k(z)| \leq \sqrt[n]{nB}$, while  \\
$|h(z)|=\sqrt[n]{\alpha^2R^2 +2\alpha \beta R \cos \theta +\beta^2}$,
leading to inequality \\
$\alpha-\beta >(nB)^{2/n}R$. \\

We collect the inequalities together, as
\[\alpha -\beta >(nB)^{2/n}R \]
\[\alpha \epsilon <\beta \]
\[(\alpha +1)\epsilon + (2B)^{1/n} <\beta \]
\[\alpha \geq (nB)^{1/n}.\]
This system can be solved if we first choose big $R$ to embrace all the roots of $f(z)$, then 
$\alpha$ to satisfy the last equation, then $\beta$, and finally $\epsilon$ as small as we need. \\

Since the linear function $h(z)$ has exactly one root in the slit domain $D$, by Rouch\'e's theorem we conclude that every function $f_k(z),\; k=0,1,\ldots , n-1,$ has exactly one root, call it 
$\xi_k$, in the complex plane. We have proved the main claim of this note.

\begin{prop} Every function $f_k(z),\; k=0,1,2,\ldots, n-1,$ has exactly one root $\xi_k$ in the complex plane. $\hfill \qed$ \end{prop}

To derive explicit representations for the roots $\xi_k$, we consider contour integrals
\beq I_k=\oint_{\Omega}\frac{z-\xi_k}{(z^2+1)f_k(z)}dz,\; k=0,1,...,n-1, \eeq
and compute each of them twice, first by the residue theorem, and then by distorting the contour of integration. For every $k$, $f_k(\xi_k)=0$, thus, the denominator of (5) has a simple pole. However, due to the numerator, at $z=\xi_k$ the integrand of (5) has a removable singularity. In addition to that, the integrand has two simple poles at $z= \pm \imath$. Computing residues at these points, we get
\[ I_k= \frac{\pi}{f_k(\imath) f_k(-\imath)}
 \left(\imath \left(f_k(\imath)+f_k(-\imath)\right) + \xi_k  \left(f_k(\imath)-f_k(-\imath)\right) \right). \]

On the other hand, letting in (5) $\epsilon \rightarrow 0$ and $R \rightarrow \infty$, the contour $\Omega$ becomes the positive $x-$axis, traversed twice, thus,
\[ I_k=\int_0^{\infty}\frac{(x-\xi_k)dx}{(x^2+1)f_k(x)} - \int_0^{\infty}\frac{(x-\xi_k)dx}{(x^2+1)f_k(xe^{2\pi \imath})} .\]

Equating the two expressions for $I_k$, we get  a linear equation for $\xi_k$. Solving this equation, we derive the explicit formulas for the roots $\xi_k,\, k=0,1,\ldots , n-1,$
{\Large
\beq \xi_k= \frac{\int_0^{\infty}\frac{x\left(f_k(xe^{2\pi \imath})-f_k(x)\right)}{(x^2+1)f_k(x)f_k(xe^{2\pi \imath})}dx - \pi \imath \frac{f_k(\imath)+f_k(-\imath)}{f_k(\imath)f_k(-\imath)}} {\int_0^{\infty}\frac{f_k(xe^{2\pi \imath})-f_k(x)}{(x^2+1)f_k(x)f_k(xe^{2\pi \imath})}dx + \pi \frac{f_k(\imath)-f_k(-\imath)}{f_k(\imath)f_k(-\imath)}}. \eeq  }

\begin{rem} These formulas are unsuitable for numerical computations, since the integrands decay very slowly. However, there are simple substitutions, which make such computations feasible even on a good calculator like TI-89, not to mention computers; see \cite{Kh1}.
\end{rem}

\bigskip

\makebox[3cm]{}      Bronx Community College \\
\makebox[3cm]{}      of The City University of New York \\
\makebox[3cm]{}      USA \\
\makebox[3cm]{}      E-mail:  alexander.kheyfits@bcc.cuny.edu \\

\end{document}